\documentclass[11pt]{amsart}

\usepackage{amsmath,amsthm,amssymb,eucal}

\newtheorem{thm}{Theorem}
\newtheorem{lem}[thm]{Lemma}
\newtheorem{prop}[thm]{Proposition}

\theoremstyle{definition}
\newtheorem{defn}{Definition}
\newtheorem{rk}{Remark}

\theoremstyle{remark}
\newtheorem{ex}{Example}

\def\epi{\mathop{\textup{epi\,}}}
\def\dist{\mathop{\textup{dist\,}}}
\def\sign{\mathop{\textup{sign\,}}}
\def\supp{\mathop{\textup{supp\,}}}
\def\divv{\mathop{\textup{div\,}}}

\def\Rset{{\mathbb R}}
\def\Nset{{\mathbb N}}

\def\dint{\displaystyle \int }
\def\dfrac#1#2{{\displaystyle \frac{#1}{#2}}}   

\numberwithin{equation}{section}

\begin{document}

\title[A nonsmooth variational approach]%
{A nonsmooth variational approach to differential  problems.
A case study of nonresonance  under the first eigenvalue
for a strongly nonlinear elliptic problem}

\author{
		Youssef JABRI
		}
\address{
     University Mohamed~I, Department of Mathematics,     
     Faculty of Sciences, Box 524, 60000 Oujda, Morocco
     }
\email{jabri@sciences.univ-oujda.ac.ma} 
\date{}

 \begin{abstract}
    We adapt a technique of nonsmooth critical point theory developed 
    by Degiovanni-Zani for a semilinear problem involving the Laplacian 
    to the the case of the $p$-Laplacian. We suppose only coercivity 
    conditions on the potential  and impose no growth condition of the 
    nonlinearity. The coercivity is obtained using similar 
    nonresonance conditions to~\cite{M-W-W} and 
    to~\cite{L-L} in two different results and using some comparison functions 
    and comparison spaces in a third one. it is also shown that 
    neither of the three theorems implies the two others.
 \end{abstract}

\maketitle
			    
\begin{center}
\end{center}


\section{Introduction}

The new nonsmooth trends in critical point theory dealing with 
continuous functionals on metric spaces 
\cite{corv-De-Ma,degiovanni-marzocchi,Katriel,ioffe-schwartzman,ioffe-schwartzman2} 
brought important developments, in particular in the abstract part of 
the theory \cite{corv1,corv2,corv3,corv4,ribarska,ribarska2}.  The 
reason is that the definition of the \emph{weak slope} they use to 
measure regularity relies on a qualitative topological property on 
level sets which avoid us the technicalities of using pseudo-gradients 
to construct deformations.  Hence, it is particularly useful to prove 
abstract critical point theorems.  While considering applications to 
boundary value problems would be, a priori, more delicate.  This is 
natural because we would need some appropriate calculus and a new 
sense for weak solutions, weaker than the traditional one, strong 
enough to be useful and still with a practical relevance.  Many early 
applications of the nonsmooth theory to several problems in partial 
differential equations and variational inequalities where the 
functionals involved are not locally Lipschitz continuous were 
given~\cite{2,4,5,6,7,13,15,19,23} but this passes in general through 
a great deal of technicality.  Recently, a specific subdifferential 
calculus has been developed by Campa and 
degiovanni~\cite{19,campa-degio} when the underlying space is normed 
and has been used successfully in applications where the Clarke's 
subdifferential is in general 
inappropriate~\cite{campa-degio,Degio-Zani2,arioli,a-g1,a-g3,squassina1,squassina2,squassina3,squassina4,squassina5}.

In the present paper, we use a general approach that combines the 
basic ideas of variational methods with these new developments and 
study nonresonance under the first eigenvalue for a strongly nonlinear 
elliptic problem involving the $p$-Laplacian operator where no growth 
condition is supposed.
The strategy followed can be described by the following scheme:
\begin{enumerate}
\item Give an acceptable definition of a weak solution to the problem 
to be treated by changing the traditional space 
$\mathcal{C}^{\infty}_{0} (\Omega)$ of test functions by an other 
dense subspace of $W^{1,p}_{0}(\Omega)$.

\item Show that a critical point is a weak solution in the sense of 
the previous point with the help of a minimum of assumptions.

\item Prove the effective existence of a critical point for the 
nonsmooth theory by one of the available abstract theorems for our 
specific nonlinear problem.

\end{enumerate}

The first two points are general and may be used with different 
boundary value problems.  We could trace their origin to 
the papers~\cite{An-Go,Degio-Zani1}.  But, they are more visible 
in~\cite{Degio-Zani2} where the authors have successfully tried a 
nonsmooth variational approach for a coercive semilinear problem with 
the Laplacian operator.
Nevertheless, these are general and may be used for different boundary 
value problems with other geometries as in the situation of linking 
results.

\section{Situating the problem}

We recall the definition of regularity in the sense of 
Degiovanni-Marzocchi~\cite{degiovanni-marzocchi} and 
Katriel~\cite{Katriel} and the key properties of the 
subdifferential calculus developed for by 
Campa-Degiovanni~\cite{campa-degio,19} we will use 
in the sequel.

Let $(E,d)$ be a metric space, $f\colon E\to\overline{\Rset}$ a 
function and denote the ball of center $u$ and radius $\sigma$ by 
$ B (u,\sigma)$ and the epigraph of $f$ by
$$
\epi(f)=\{(u,\lambda)\in E\times\Rset;\ f(u)\leq \lambda\}.
$$
Consider $\epi(f)$ as a subspace of the metric space $E\times\Rset$ 
for the metric $\dist((u,\lambda),(u',\lambda')) =\big(\dist 
(u,u')^{2} +(\lambda'-\lambda)^{2} \big)^{1/2}$.
\begin{defn}
	 A point $u\in E$ such that $f(u)\in\Rset$ is
	\emph{${\delta}$-regular} for a $\delta>0$, if there exists 
	$\sigma>0$ and a continuous deformation $\nu\colon 
	\big( B_{\sigma}(u,f(u)) \cap \epi(f)\big) \times [0,\sigma]\to E$ such that for all 
	$(w,\lambda)\in B(u,f(u)) \cap \epi(f)$ and all $t\in[0,\sigma]$: 
	$$
	\dist(\nu((w,\mu),t),w)\leq t\qquad\hbox{and}\qquad 
	f(\nu((w,t),t))\leq \mu-\delta t.
	$$
	The point $u$ is \emph{regular} if there exists $\delta>0$ such that $u$ is
	$\delta$-regular and it is  \emph{critical} if it is not regular.
\end{defn}
When $f$ is continuous, this definition may be simplified to the 
following form.
\begin{defn}
    Let $f\colon E\to{\Rset}$ be a continuous function.  A point $u\in E$ 
    is \emph{${\delta}$-regular} if there exists 
    a continuous deformation $\eta\colon [0,\sigma] \times  B(u,\sigma) 
    \to E$ such that for all $(t,u)\in[0,\sigma]\times  B(u,\sigma)$:
	$$
	\dist(\eta(t,u),u)\leq t\qquad\hbox{and}\qquad 
	f(u)-f(\eta(t,u))\geq \delta t.
	$$
\end{defn}

The regularity of $u$ is evaluated by a generalized notion of the norm 
of the derivative called by Degiovanni-Marzocchi and Katriel 
respectively \emph{weak slope} and \emph{regularity constant} of $f$ 
at~$u$:
$$
|df|(u)=\sup\big\{ \delta>0;\ u\text{ is }\delta-\text{regular} 
\big\},
$$
which reduces to $||f'(u)||$ when $f$ is of class $\mathcal{C}^{1}$.

When $u$ is a local minimum of $f$, it is a critical point, i.e., 
$\vert df\vert(u)=0$.  And when $f$ satisfies the geometric conditions 
of the mountain pass theorem of Ambrosetti-Rabinowitz 
\cite{Ambro-Rabin}, and an appropriate form of the Palais-Smale 
condition, it has a critical point whose value is characterized by the 
usual $\inf\max$ argument (cf.  \cite{degiovanni-marzocchi,Katriel}).

When $X$ is a normed space.

\begin{defn}
For any $u\in X$ such that $f(u)\in\Rset$, $v\in X$ and 
$\varepsilon>0$, we define $f^{0}_{\varepsilon}(u;v)$ as the infimum 
of $r\in\Rset$ such that there exists $\delta>0$ and a continuous 
mapping
$$
\nu\colon \big( B_{\delta}(u,f(u))\cap \epi(f)\big)\times ]0,\delta]\to 
 B_{\varepsilon}(v)
$$            
satisfying
$$
f(w+t\nu((w,\mu),t))\leq \mu+rt
$$
for all $(w,\mu)\in  B_{\delta}(u,f(u))\cap \epi(f)$ and 
$t\in]0,\delta]$.

Define also
$$
f^{0}(u;v)=\sup_{\varepsilon>0}f^{0}_{\varepsilon}(u;v).
$$
\end{defn}
The function $f^{0}(u;.)$ is convex, lower semicontinuous (l.s.c.) 
and positively homogeneous of degree 1 (cf.  \cite{campa-degio}).

\begin{defn}
For all $u\in X$ such that $f(x)\in\Rset$, we set
$$
\partial f(u)=\big\{ \alpha\in X^{*};\ \langle \alpha,v\rangle \leq 
f^{0}(u;v)\ \forall v\in X\big\}.
$$
\end{defn}
\noindent The real $f^{0}(u;v)$ is larger than the 
Clarke-Rockafellar generalized directional derivative.   
Hence, $\partial f(u)$ contains the Clarke's subdifferential 
of $f$ at $u$ which may be empty when $\partial f(u)$ is not, as we 
can see for example for $f(u)=u-\sqrt{|u|}$ at $u=0$ (there $|df|(0)=0$).  
Nevertheless they agree when $f$ is locally Lipschitz.
The notions $f^{0}(u;v)$ and
$\partial f(u)$ have been introduced in \cite{campa-degio,19} and
are more adapted to use with the weak slope of Degiovanni-Marzocchi.

\begin{prop}\label{prop1}
If $u\in X$ is such that $f(u)\in\Rset$, then\\
\textup{(a)} $|df|(u)<+\infty \iff \partial f(u)\not= \varnothing$.\\
\textup{(b)} $|df|(u)<+\infty \Rightarrow |df|(u)\geq \min\{ ||\alpha||;\ 
\alpha\in\partial f(u)\}$.
\end{prop}

\begin{rk}
	 Examples where (b) is strict exist.
\end{rk}

\begin{prop}\label{prop2}
When $u\in X$ is such that $f(u)\in\Rset$ and $g\colon X\to\Rset$ is 
Lipschitz continuous, then\\  
$$
\partial (f+g)(u)\subset \partial (f)(u)+\partial (g)(u).
$$
If moreover, $g$ is of class $\mathcal{C}^{1}$, we have equality 
and $\partial g(u)=\{g'(u)\}$.
\end{prop}

To illustrate the abstract approach, we will study nonresonance under 
the first eigenvalue for a strongly nonlinear elliptic problem for the 
$p$-Laplacian operator as in~\cite{Ch}.  The former reference contains 
improvements to some results in~\cite{An-Go} adapting the contents 
of~\cite{Degio-Zani1} to the case of the $p$-Laplacian.  Both results 
might be adapted to a nonsmooth variational context.  Degiovanni-Zani 
adapted their results of~\cite{Degio-Zani1} in \cite{Degio-Zani2}.  We 
will do the same, by analogy to \cite{Degio-Zani2} to the results 
of~\cite{Ch} extending \cite{Degio-Zani2} to the case where $p$ may be 
different from 2.

The usual norm of $W^{1,p}_{0}(\Omega)$ will be denoted $||.||_{1,p}$, 
the Lebesgue measure of a set $A$ is $|A|$ while 
$ B(x,\rho)= B_{\rho}(x)$ is the open ball of center $x$ and 
radius $\rho$.  The notation $m'$ stands for the conjugate of $m$, 
i.e.  $1/m + 1/m' =1$.

The  problem to be treated is:
$$
(\mathcal{P})~\left\{
\begin{array}{rll}
-\Delta_{p} u = & f(x,u) +h & \text{in }\Omega\\
u=& 0& \text{on }\partial\Omega
\end{array}
\right.
$$
where $\Omega$ is a bounded domain of $\Rset^N$, $\Delta_{p}\colon 
W^{1,p}_{0}(\Omega)\to W^{-1,p'}(\Omega)$ is the 
$p$-Laplacian operator defined by 
$$
 -\Delta_{p}u \equiv \divv
\left(|\nabla u|^{p-2}\nabla u\right), \qquad 1<p<\infty.
$$
The $p$-Laplacian is a degenerated quasilinear elliptic operator that 
reduces to the classical Laplacian if $p=2$.
The notation $\langle .,.\rangle$ stands hereafter for the duality pairing 
between $ W^{-1,p'}(\Omega)$ and $ W^{1,p}_{0}(\Omega)$.
While $f\colon \Omega\times \Rset \to \Rset$ is a Carath\'eodory 
function and $h\in W^{-1,p'}(\Omega)$.

Consider the energy functional $\Phi\colon W^{1,p}_0(\Omega) \to 
\overline{\Rset}$ associated to the problem
$$
\Phi(u)=\frac{1}{p}\int_\Omega |\nabla u |^p\,dx-\int_\Omega 
F(x,u)\,dx -\langle h,u\rangle,
$$
where $F(x,s)=\int_0^s f(x,t)\,dt$. We are interested in conditions to 
be imposed on the nonlinearity $f$ in order that problem $(\mathcal{P})$ 
admits at least one solution $u(x)$ for any given $h(x)$. Such 
conditions are usually called \emph{nonresonance conditions}.

When the nonlinearity satisfies a growth condition of the type:
\begin{equation}\label{(f)}
|f(x,s)|\leq a |s|^{q-1} +b(x) \text{ for all }s \in \Rset, 
\text{ and a.e. in }\Omega, 
\end{equation}
with $q<p^*$ where the Sobolev exponent $p^{*}=\frac{Np}{N-p}$ when 
$p<N$ and $p^{*}=+\infty$ when $p\geq N$ and $b(x)\in L^{(p^{*})'}(\Omega)$, 
the functional $\Phi$ is well defined, of class $\mathcal{C}^1$, 
l.s.c.  and its critical points are weak solutions of $(\mathcal{P})$ 
in the usual sense.

But when this growth condition is not satisfied, $\Phi$ is not 
necessarily of class $\mathcal{C}^1$ on $W^{1,p}_0(\Omega)$ and may 
take infinite values.

The first eigenvalue of the $p$-Laplacian  characterized by the 
variational formulation
$$
\lambda_{1}=\lambda_{1}(-\Delta_{p})=\min \left\{ 
\frac{\int_{\Omega}|\nabla u|^{p}}{\int_{\Omega}|u|^{p}}\,dx;\ u\in 
W^{1,p}_{0}(\Omega)\setminus \{ 0\}\right\}
$$
is known (cf.  \cite{An} for example) to be associated to a simple 
eigenfunction that does not change sign.  We will denote by 
$\varphi_{1}$ the normalized eigenfunction such that $\varphi_{1}>0$ 
almost everywhere.

A procedure used to treat $(\mathcal{P})$ when the nonlinearity lies 
asymptotically on the left of $\lambda_{1}$ consists in supposing a 
``coercivity'' condition on $F$ of the type:
\begin{equation}\label{(F)}
\limsup_{s\to\pm\infty}\frac{pF(x,s)}{|s|^p}<\lambda_1\qquad 
\text{for almost every }x\in\Omega
\end{equation}
and minimizing $\Phi$ on $W^{1,p}_{0}(\Omega)$.  The minimum being a 
weak solution of $(\mathcal{P})$ in an appropriate sense (cf.  
\cite{An-Go,Degio-Zani1,Degio-Zani2,Ch}).

An other way, is to obtain \emph{a priori} estimates on the solutions 
of some equations approximating $(\mathcal{P})$ and to show that their 
weak limit is indeed a weak solution.

Notice that with the help of the conditions \eqref{(f)} and \eqref{(F)}, we know 
since the work of Hammerstein (1930) that $(\mathcal{P})$ admits a 
weak solution that minimizes the functional $\Phi$ on 
$W^{1,p}_{0}(\Omega)$.

The condition \eqref{(F)} does not imply a growth condition on $f$ unless 
$f(x,u)$ is convex in $u$ (see for example 
\cite{M-W1,jabri-moussaoui2}).

In~\cite{An-Go}, Anane and Gossez supposed only a one-sided growth 
condition with respect to the Sobolev (conjugate) exponent that do not 
suffice to guarantee the differentiability of $\Phi$, which may even 
take infinite values.
Nevertheless, they showed that any minimum of $\Phi$ solves 
$(\mathcal{P})$ in a suitable sense.

Then, Degiovanni-Zani~\cite{Degio-Zani1} in the case $p=2$ and 
Chakrone in~\cite{Ch} for $1<p<\infty$ supposed only that $f$ maps 
$L^\infty(\Omega)$ to $L^1(\Omega)$:
\begin{equation}\label{(f_0)}
\sup_{|s|\leq R}|f(.,s)|\in L^1_{loc}(\Omega),\qquad \forall R>0
\end{equation}
and a coercivity condition of the type \eqref{(F)}.  They proved that 
any minimum $u$ of $\Phi$, which is not of class $\mathcal{C}^1$ on 
$W^{1,p}_{0}(\Omega)$ and may take infinite values too, is a weak 
solution of $(\mathcal{P})$ in the sense
$$
\int_\Omega |\nabla u|^{p-2}\nabla u \nabla v\,dx =\int_\Omega 
f(x,u)v\,dx +\langle h,v\rangle,
$$
for $v$ in a dense subspace of $W^{1,p}_0(\Omega)$.
Recently, in 1998 Degiovanni and Zani adapted their technique of 
\cite{Degio-Zani1} for the Laplacian operator ($p=2$) to a nonsmooth 
analytical context in \cite{Degio-Zani2}.  We will do here the same 
with the results of \cite{Ch} for the $p$-Laplacian ($1<p<\infty$).

In the autonomous case $f(x,s)=f(s)$, De~Figueiredo and 
Gossez~\cite{Defi-Go} have proved the existence of solutions for any $h\in 
L^{\infty}(\Omega)$ by a topological method. They supposed only a 
coercivity condition and established that
$$
\int_{\Omega}|\nabla u|^{p-2}\nabla u\nabla v\,dx=\int_{\Omega}f(x,u) 
v\,dx +\langle h,v\rangle
$$
for all $v\in W^{1,p}_{0}(\Omega)\cap L^{\infty}(\Omega) \cup \{u\}$ 
but the solution obtained may not minimize $\Phi$.  Indeed, an example 
is given in~\cite{Defi-Go} in the case $p=2$ and an other one is given 
in~\cite{Ch} where $p$ may be different from~2.

Notice that in our case, the condition \eqref{(f_0)} implies no growth 
condition on $f$ as it may be seen in the following example 
from~\cite{Ch}.
\begin{ex}
	Consider the function 
	$$
	f(x,s)=\left\{
	\begin{array}{ll}
	d(x)\left( \sin \left( 
	\dfrac{\pi s}{2}\right)-\dfrac{\sign(s)}{2}\right) \exp \left( 
	\dfrac{2\cos\left( \frac{\pi s}{2}\right)}{\pi} +\dfrac{|s|-1}{2} 
	\right) & \text{if } |s|\geq 1\\
	d(x)\dfrac{s}{2}(10 s^{2}-9) &\text{if } |s|\leq 1
	\end{array}
	\right.
	$$
	where $d(x)\in L^{1}_{loc}(\Omega)$ and $d(x)\geq 0$ almost 
	everywhere in $\Omega$, so that
	$$
	F(x,s)=\left\{
	\begin{array}{ll}
	-d(x)\exp\left(\dfrac{2\cos \left( 
	\frac{\pi s}{2}\right)}{\pi} \right) \exp \left( 
	\dfrac{|s|-1}{2} 
	\right) & \text{if } |s|\geq 1\\
	-d(x)\dfrac{s^{2}}{4}(-5 s^{2}+9) &\text{if } |s|\leq 1
	\end{array}
	\right. 
	$$
	Then, $F(x,s)\leq 0$ for all $s\in\Rset$ almost everywhere in 
	$\Omega$. So, $\Phi$ is coercive. Nevertheless, as we can check 
	easily, $f$ satisfies no growth condition.
\end{ex}

\section{Theoretical approach}

We will show that when \eqref{(f_0)} is fulfilled, a critical point in the 
sense of Degiovanni-Marzocchi~\cite{degiovanni-marzocchi}, is a weak 
solution of $(\mathcal{P})$ in an acceptable sense.

\begin{defn}
For $1\leq p\leq \infty$, the space $L^p_0(\Omega)$ is defined 
by:
$$
L^p_0=\big\{v\in L^{p}_0;\ v(x)=0\text{ a.e.  outside a compact 
subset of }\Omega \big\}.
$$

Consider $u\in W^{1,p}_0(\Omega)$, we set
$$
V_{u}=\big\{v\in W^{1,p}_{loc};\ u\in L^\infty(\{x\in 
\Omega;\ v(x)\not=0\}) \big\}.
$$
\end{defn}

\begin{prop}[Brezis-Browder {\cite{B-B}}]\label{bre-bro-prop}
If $u\in W^{1,p}_0(\Omega)$, there exists a sequence $(u_n)_n\subset 
W^{1,p}_0(\Omega)$ such that:\\  
\textup{(i)} $(u_n)_n\subset W^{1,p}_{0}(\Omega)\cap L^{\infty}_{0}(\Omega)$.\\
\textup{(ii)} $|u_n(x)|\leq |u(x)|$ and $u_n(x).u(x) \geq 0$ a.e. in $\Omega$.\\
\textup{(iii)} $u_n\to_{n\to\infty}u$ in $W^{1,p}_0(\Omega)$.
\end{prop}

The linear space $V_{u}$ enjoys some nice properties.
The next proposition that refines in some sense the former one is due 
to~\cite{Degio-Zani1} when $p=2$.
\begin{prop}[Chakrone \mbox{\cite{Ch}}]
The space $V_u$ is dense in $W^{1,p}_0(\Omega)$.  And if we suppose 
suppose that \eqref{(f_0)} holds, then
$$
A_u=\big\{\varphi\in W^{1,p}_0(\Omega);\ f(x,u)\varphi\in L^1(\Omega) \big\}
$$
is a dense subspace of $W^{1,p}_0(\Omega)$ as $V_u\subset A_u$.  More 
precisely, Brezis-Browder's result holds true if we replace 
$W^{1,p}_{0}(\Omega)\cap L^{\infty}_{0}(\Omega)$ by $A_{u}$.
\end{prop}
\begin{proof}
	It suffices to show that $V_{u}$ is dense in $W^{1,p}_{0}(\Omega)$ 
	and that $V_{u}\subset A_{u}$ when \eqref{(f_0)} holds.\\
	\textbf{1. The density of $V_{u}$ in $W^{1,p}_{0}(\Omega)$:}\\
	 We have to show that for any $\varphi\in W^{1,p}_{0}(\Omega)$, there 
	 exists a sequence $(\varphi_{n})_{n}\subset V_{u}$ satisfying 
	(ii) and (iii). This is done in two steps. First, we show it is true 
	for all $\varphi\in W^{1,p}_{0}(\Omega)\cap L_{0}^{\infty}(\Omega)$.
	Then, using Proposition~\ref{bre-bro-prop}, we show it is true in 
	$W^{1,p}_{0}(\Omega)$.
	
	\textbf{First Step:} Suppose $\varphi\in W^{1,p}_{0}(\Omega)\cap 
	L_{0}^{\infty}(\Omega)$ and consider a sequence 
	$(\Theta_{n})_{n}\subset \mathcal{C}^{\infty}_{0}(\Rset)$ such that:\\
	(1) $\supp \Theta_{n}\subset [-n,n]$,\\
	(2) $\Theta\equiv 1$ on $[-n+1,n-1]$,\\
	(3) $0\leq \theta_{n}\leq 1$ on $\Rset$ and\\
	(4) $|\Theta_{n}'(s)|\leq 2$.
	
	The sequence we are looking for is obtained by setting
	$$
	\varphi_{n}(x)=(\Theta\circ u)(x)\varphi(x) \quad\text{for a.e. }x \text{ 
	in }\Omega.
	$$
	Indeed, let's check the following three points\\
	(a) $\varphi_{n}\in V_{u}$,\\
	(b) $|\varphi_{n}(x)|\leq |\varphi(x)|$ and 
	$\varphi_{n}(x)\varphi(x)>0$ 
	a.e. in $\Omega$ and\\
	(c) $\varphi_{n}\to \varphi$ in $W^{1,p}_{0}(\Omega)$.\\
	For (a), since $\varphi\in L^{\infty}_{0}(\Omega)$, we have that 
	$\varphi_{n}\in L^{\infty}_{0}(\Omega)$ and it's clear by~(4) that 
	$\varphi_{n}\in W^{1,p}_{0}(\Omega)$. Finally, by (1), 
	$u(x)\in[-n,n]$ for a.e. $x$ in $\{ x\in \Omega;\ 
	\varphi_{n}(x)\not= 0\}$.
	
	The assumption (b) is a consequence of (3).
	
	For (c), by (2), $\varphi_{n}(x)\to\varphi(x)$ a.e. in $\Omega$ and
	$$
	\frac{\partial \varphi_{n}}{\partial x_{i}}(x)=\Theta'_{n}(u(x))\frac{\partial 
	u}{x_{i}}\varphi(x)+\Theta_{n}(u(x))\frac{\partial \varphi}{\partial 
	x_{i}}\ \to\ \frac{\partial \varphi}{\partial x_{i}}\text{ in }\Omega.
	$$
	And by (4), 
	$$
	\left|\frac{\partial \varphi_{n}}{\partial x_{i}}(x) 
	\right|\leq 
	\left|\frac{\partial u}{\partial x_{i}}(x)\right|  |\varphi(x)|
	+\left|\frac{\partial \varphi}{\partial x_{i}}(x) 
	\right|\in L^{p}(\Omega).
	$$
	
	Finally, by the dominated convergence theorem we get (c).
	
	\textbf{Second Step:} Suppose that $\varphi\in W^{1,p}_{0}(\Omega)$. By 
	Proposition~\ref{bre-bro-prop}, there is a sequence 
	$(\psi_{n})_{n}\subset W^{1,p}_{0}(\Omega)$ satisfying (i), (ii) 
	and (iii).
	
	For $k=1,2,\ldots$, there is $n_{k}\in \Nset$ such that 
	$||\psi_{n_{k}}-\varphi||_{1,p}\leq {1}/{k}$. Since 
	$\psi_{n_{k}}\in W^{1,p}_{0}(\Omega)\cap 
	L^{\infty}_{0}(\Omega)$, by the first step, there is $\psi_{k}\in 
	V_{u}$ such that $|\varphi_{k}(x) |\leq |\psi_{n_{k}}(x)|$ and 
	$\varphi_{k}(x)\psi_{n_{k}}(x)\geq 0$ almost everywhere in $\Omega$ 
	and $||\varphi_{k}-\psi_{n_{k}}||_{1,p}\leq {1}/{k}$, so that 
	$(\varphi_{k})_{k}$ is the sequence we are seeking. Indeed, 
	$|\varphi_{k}(x) |\leq |\psi_{n_{k}}(x)|\leq |\varphi(x) |$, 
	$\varphi_{k}(x)\varphi(x)\geq 0$ a.e. in $\Omega$ and 
	$||\varphi_{k}-\varphi(x)||_{1,p}\leq ||\varphi_{k}-\psi_{n_{k}}||_{1,p}
	+ ||\psi_{n_{k}}-\varphi(x)||_{1,p}\leq {2}/{k}$.\\[2mm]
	\textbf{The inclusion $V_{u}\subset A_{u}$:}\\
	Indeed, for $\varphi\in V_{u}$, set $E=\big\{ x\in \Omega;\ 
	\varphi(x)\neq 0\big\}$ so that 
	$$
	\begin{array}{rl}
	|f(x,u)\varphi | & = \left| f(x,u) \chi_{E}\varphi(x)\right|\\[2mm]
	              & \leq \max \big\{ | f(x,s) \varphi(x)|; |s|\leq||u 
	              ||_{L^{\infty}(E)}\big\}
	\end{array}
	$$
	where $\chi_{E}$ is the characteristic function of the set $E$.
	
	By \eqref{(f_0)}, the last term lies to $L^{1}(\Omega)$, so that 
	$\varphi\in A_{u}$.
\end{proof}

So, we have that
\begin{equation}
\begin{split}
\Lambda u=\sup\left\{\int_\Omega f(x,u)v\,dx;\ v\in A_u, ||v||_{1,p}\leq 
1 \right\}=\\
\sup\left\{\int_\Omega f(x,u)v\,dx;\ v\in V_u, ||v||_{1,p}\leq 1 
\right\}.
\end{split}
\end{equation}

\begin{defn}\label{defin}
Consider $u\in W^{1,p}_{loc}(\Omega)$, we say that $f(x,u)\in 
W^{-1,p'}(\Omega)$ if 
$$
\Lambda u<+\infty.
$$
\end{defn}

Then, the mapping $T\colon V_u\to\Rset$ defined by
$$
T(\varphi)= \int_\Omega f(x,u)\varphi\,dx
$$ 
is linear, continuous and admits an extension $\tilde{T}$ to the whole 
space $W^{1,p}_0(\Omega)$.  Henceforth, we will make the identification 
$f(x,u)=\tilde{T}$; this justifies the terminology of 
Definition~\ref{defin}.

\begin{defn}[Weak solution]
A point $u\in W^{1,p}_0(\Omega)$ is a weak solution of $(\mathcal{P})$ 
if  $f(x,u)\in W^{-1,p'}(\Omega)$ (in the sense of 
Definition~\ref{defin}) and $(\mathcal{P})$ is satisfied in 
$W^{-1,p'}(\Omega)$.

In particular, we would have
$$
\int_\Omega |\nabla u|^{p-2}\nabla u \nabla v\,dx =\int_\Omega 
f(x,u)v\,dx + \langle h,v\rangle,\quad \forall v\in A_u.  
$$
\end{defn}
According to the conventions of \cite{rockafellar}, adopted also in 
\cite{Degio-Zani2}, we will consider that by definition
$$
\int_{\Omega} F(x,u)^{+}\,dx = \int_{\Omega} F(x,u)^{-}\,dx =+\infty 
\quad\text{implies that }\int_{\Omega} F(x,u)\,dx =\infty.
$$
And since the definition of the functional $\Phi$ contains the term 
$\int_{\Omega} F(x,u)\,dx$, the above convention has to be considered 
for $-F$.
\begin{thm}\label{thmfond1}
	Let $u\in W^{1,p}_{0}(\Omega)$ such that
	$\Psi(u)=\int_{\Omega}F(x,u)\,dx \in\Rset$, then\\
	\textup{\textbf{(a)}} For any $v\in W^{1,p}_{0}(\Omega)$ such that $f(x,u)v\in 
	L^{1}(\Omega)$, we have
	$$
	\Psi^{\circ}(u;v) \leq \int_{\Omega} f(x,u)v\,dx.
	$$
	This is true in particular for all $v\in V_{u}$.\\
	\textup{\textbf{(b)}} If $\partial \Psi(u)\not=\varnothing$, 
	then $f(x,u)\in 
	W^{-1,p'}(\Omega)$ and $\partial \Psi(u)=\{ f(x,u)\}$.
\end{thm}
\begin{proof}
\textbf{(a)} 
We will prove that the assertion holds for all $v\in V_{u}$ in a first step and 
then, for all $v\in A_u$ in a second step.
	
\noindent  \textbf{1 $^{\text{st}}$ step:}\\
	Consider $v\in V_{u}$, $\varepsilon>0$, $R>0$ and
	$$
	r>\int_{\Omega}f(x,u)v\,dx.
	$$
	Consider also a smooth function $\Theta_{R}\colon \Rset \to [0,1]$  
	with support in $[-2R,2R]$ such that $\Theta_{R}(s)=1$ 
	in $[-R,R]$ and $|\Theta'_{R}(s)| \leq 2/R$ in $\Rset$. For $R$ large
	enough, we have $||\Theta_{R}(u)v-v||_{1,p}<\varepsilon$ and
	$$
	r>\int_{\Omega}f(x,u)\Theta(u) v\,dx.
	$$
	\textbf{Affirmation:}
	$$
	\lim_{\substack{
	w\to u\\
	t\to 0^{+}
	}
	}\int_{\Omega}\frac{F(x,w+t\Theta_{R}(w)v)-F(x,w)}{t}\,dx 
	=\int_{\Omega}f(x,u)\Theta_{R}(u)v\,dx.
	$$
	Indeed, consider $(w_{h})_{h}$ such that $w_{h}\to u$ in $W^{1,p}_{0}(\Omega)$ 
	and $(t_{h})_{h}$ such that $t_{h}\to 0^{+}$. Without loss of
	generality, we may suppose that $w_{h}(x)\to u(x)$ almost
	everywhere and that $t_{h}\leq 1$.
	It follows that
	$$
	\lim_{h}\frac{F(x,w_{h}+t_{h}\Theta_{R}(w_{h})v)-F(x,w_{h})}{t_{h}} 
	=f(x,u)\Theta_{R}(u)v \qquad\text{ a.e. in }\Omega.
	$$
	On the other hand, for almost every $x$ in $\Omega$, there exists 
	$\tau_{h}\in ]0,t_{h}[$ such that
	$$
	\begin{array}{ll}
		\dfrac{F(x,w_{h}+t_{h}\Theta_{R}(w_{h})v)-F(x,w_{h})}{t_{h}}  & 
		=|f(x,w_{h}+t_{h}\Theta_{R}(w_{h})v)|\;|v|\Theta_{R}(w_{h})  \\
		 & \leq \left(\sup_{|s|\leq T}|f(x,s)| \right) |v|
	\end{array}
	$$
	where $T=2R+||v||_{\infty}$. The affirmation follows then from~\eqref{(f_0)} 
	and Lebesgue's dominated convergence theorem.
	
	Therefore, there exists $\delta>0$ such that for any $w\in 
	W^{1,p}_{0}(\Omega)$ satisfying $||w-v||_{1,p}<\delta$ and 
	$0<t<\delta$ we have
	$$
		||\Theta_{R}(w)v-v||_{1,p}< \varepsilon, 
	$$
	$$
		F(x,w+t\Theta_{R}(w)v)-F(x,w) \in L^{1}(\Omega),  
	$$
	and
	$$
		\int_{\Omega}\frac{F(x,w+t\Theta_{R}(w)v)-F(x,w)}{t}\,dx <r.
	$$
	Consider now the  continuous function
	$
	\nu\colon 
	( B_{\delta}(u,\Psi(u))\cap \epi(\Psi)) \times 
	]0,\delta] \to  B_{\varepsilon}(v)
	$ 
	defined by 
	$\nu((w,\mu),t)=\Theta_{R}(w)v$.
	
	We have that $\Psi(w+t\nu((w,\mu),t))=-\infty$ if and only 
	if  $\Psi(w)=-\infty$ and $\Psi(w+t\nu((w,\mu),t))\leq 
	\Psi(w)+rt\leq \mu+rt$. 
	
	Therefore,
	$\Psi^{\circ}_{\varepsilon}(u;v)\leq r$. Since $\varepsilon$ 
	has been taken arbitrary, we have
	$\Psi^{\circ}(u;v)\leq r$.  
	
	And since $r\in\Rset$ has been also taken arbitrary, we get \textbf{(a)} 
	for all $v\in 
	V_{u}$.  

	\bigskip

	\noindent \textbf{2 $^{\text{nd}}$ step:}\\
	Let $v\in W^{1,p}_{0}$ such that $f(x,u)v\in L^{1}(\Omega)$, i.e $u\in 
	A_{u}$. By Proposition~\ref{bre-bro-prop}, there exists 
	$(v_{n})_{n}\subset V_{u}$ such that $v_{n}\to v$ in 
	$W^{1,p}_{0}(\Omega)$. 
	We can suppose that $v_{n}\to v$ almost everywhere, 
	$$
	|f(x,u)v_{n}|\leq |f(x,u)v|\text{ and }f(x,u)v_{n}\to 
	f(x,u)v\qquad\text{ a.e. }. 
	$$
	By the dominated convergence theorem, we have then
	$$
	\lim_{h}\int_{\Omega}f(x,u)v_{n}\,dx =\int_{\Omega}f(x,u)v\,dx.
	$$
	And since $\Psi^{\circ}(u;.)$ is l.s.c., the assertion \textbf{(a)} holds 
	true for all $v\in A_u$.
	
	\noindent\textbf{(b)} Let $\alpha\in \partial \Psi(u)$. Suppose that 
	$v\in W^{1,p}_{0}(\Omega)$ and that $f(x,u)v \in L^{1}(\Omega)$, then 
	$$
	\langle \alpha,v\rangle \leq \Psi^{\circ}(u;v)\leq \int_{\Omega}f(x,u)v\,dx.
	$$
	As we can change in the above inequality $v$ by $-v$, it follows that $\langle 
	\alpha,v\rangle=\dint_{\Omega}f(x,u)v\,dx$, so $f(x,u)\in 
	W^{-1,p'}(\Omega)$ and $\alpha=f(x,u)$.
\end{proof}
	
	So, we get as an immediate consequence of Propositions \ref{prop1} and~\ref{prop2} 
	and Theorem~\ref{thmfond1}, the following result.

\begin{thm}
	If $u\in W^{1,p}_{0}(\Omega)$ is a critical point of $\Phi$ in the 
	sense 
	of Degiovanni-Marzocchi such that $\Psi(u)\in\Rset$, then $u$ 
	is a weak solution of the problem $(\mathcal{P})$.
\end{thm}

\section{Coercive problems}

We will see now in an analogous approach to \cite{Degio-Zani2} some 
conditions due to~\cite{Ch} that guarantee the existence of a global 
minimum $u$ of $\Phi$ in $W^{1,p}_{0}(\Omega)$ and such that 
$\Psi(u)\in\Rset$ and hence are exactly what we need.  Three results 
are obtained using similar nonresonance conditions respectively to 
Mawhin-Ward-Willem~\cite{M-W-W} and to Landesman-Lazer~\cite{L-L} in 
the two first ones and using some comparison functions and comparison 
spaces in the third.  The three theorems are shown to be 
incomparable.

Set $G(x,s)=F(x,s)-\lambda_{1}{|s|^{p}}/{p}$.  The conditions will 
port on $G(x,s)/|s|^{\alpha}$ for $1\leq \alpha\leq p$.  Anane and 
Gossez introduced in~\cite{An-Go} some comparison functions and 
comparison spaces. Other ones are used here.
\begin{defn}
	A  continuous even function $\varphi\colon \Rset\to \Rset^{+}$ is 
	called a  comparison function  of order $\alpha$, where $1\leq \alpha\leq p$ 
	if\\
	(i) $\dfrac{\varphi(s)}{|s|^{p}}\to 0$ when $s\to +\infty$.\\
	(ii) $\dfrac{\varphi(s)}{|s|}\to +\infty$ when $s\to +\infty$.\\
	(iii) $\dfrac{\varphi(s)}{\varphi(t_{n})}\to r^{\alpha}$ when $\dfrac{s_{n}}{t_{n}}\to 
	r>0$, $s_{n}\to \infty$ and $t_{n}\to \infty$.\\
	(iv) For all $\beta>\alpha$, there exist $t_{0},a$ and $b$ such that 
	$$
	\frac{\varphi(ts)}{\varphi(t)}\leq a |s|^{\beta}+b\qquad \text{forall } t\geq t_{0}
	\text{ and all }s\geq 0.
	$$
\end{defn}
\begin{ex}~\\
	a) $\varphi(s)= |s|^{\alpha}$, $1<\alpha<p$.\\
	b) $\varphi(s)= \dfrac{|s|^{\alpha}}{\big|\log|s|\big|}$,  $1<\alpha<p$.\\
	c) $\varphi(s)= |s|^{\alpha}\big|\log|s|\big|$,  $1\leq\alpha<p$.\\
\end{ex}

\begin{defn}
	Let $1\leq\alpha\leq p$.  We denote by $X_{\alpha}$ 
	the set of all  measurable functions $\eta(x)$ on $\Omega$ 
	satisfying:\\
	(i) $\eta(x)\in L^{1}(\Omega)$  if $p=N$.\\
	(ii) $\eta(x)\in L^{q}(\Omega)$ for some $q>1$ if $p=N$.\\
	(iii) $\eta(x)\in L^{q}(\Omega)$ for some $q>(p^*/\alpha)'$  if $p<N$.\\
	The set $Y_{\alpha}$ is defined the same way as  $X_{\alpha}$ except that it 
	is required that
	 $\eta(x)\in L^{(p^*/\alpha)'}(\Omega)$ if $p<N$.
\end{defn}
For a comparison function $\varphi$ of order $\alpha$, $1\leq 
\alpha\leq p$, we denote by
$$
G^{\pm}_{\varphi}(x)=\limsup_{s\to\pm\infty} 
\frac{G(x,s)}{\varphi(s)}\text{ for almost every }x\in \Omega.
$$ 
If $\varphi(s)=|s|^{\alpha}$, we write only
$$
G^{\pm}_{\varphi}(x)=G^{\pm}_{\alpha}(x)
$$
Consider now a slightly stronger condition than \eqref{(f_0)},
\begin{equation*}
	\sup_{|s|\leq R}|f(.,s)|\in L^1(\Omega),\qquad \forall R>0
	\tag{$f_{0}$}
\end{equation*}
We have then the following coercivity results for the problem $(\mathcal{P})$:
\begin{thm}\label{thm:7}
	Suppose $(f_0)$,\\ 
	$
	{(G_{1})}\qquad
	G^{\pm}_{\alpha}\leq 0 \text{ a.e. in }\Omega
	$, \\
	$
	{(G_{1}')}\qquad
	|\{ x\in \Omega;\ G^{+}_{p}< 0 \}|\not=0\text{ and } |\{ x\in \Omega;\ 
	G^{-}_{p}< 0 \}|\not=0
	$.\\
Then, $\Phi$ achieves its minimum in a point $u$ in $W^{1,p}_{0}(\Omega)$ and
$\Psi(u)\in\Rset$. And hence $u$ is a weak solution of~$(\mathcal{P})$.
\end{thm} 

 \begin{thm}\label{thm:3}
	Suppose $(f_0)$, \\
	$
	{(G_{3})}\qquad
	G^{\pm}_{1}\leq \eta \text{ a.e. uniformly in 
	}\Omega \text{ for some } \eta\in Y_{1}
	$,
	and\\
	$
	{(G_{3}')}\qquad
	\int_{\Omega}G^{-}_{1}(x)\varphi_{1}(x)< \langle 
	h,\varphi_{1}\rangle < -\int_{\Omega}G^{+}_{1}(x)\varphi_{1}(x)$.\\
Then, the conclusion of Theorem~\ref{thm:7} holds. 
\end{thm}

\begin{thm}\label{thm:second}
	Suppose that $\varphi$ is a  comparison function of order $\alpha$, $1\leq 
	\alpha\leq p$. Suppose also that $(f_0)$ holds and  \\
	$
	{(G_{2})}\qquad
	G^{\pm}_{\varphi}\leq \eta \text{ a.e. uniformly in 
	}\Omega \text{ for some } \eta\in X_{\alpha}
	$,
	and\\
	$
	{(G_{2}')}\qquad
	\int_{\Omega}G^{+}_{\varphi}(x)(\varphi_{1}(x))^{\alpha}<0\text{ and }
	\int_{\Omega}G^{-}_{\varphi}(x)(\varphi_{1}(x))^{\alpha}<0
	$.\\
Then, the conclusion of Theorem~\ref{thm:7} holds.
\end{thm}

The proofs of the three theorems use the same technique. To give the 
best idea on the role played by comparison spaces and functions, we 
will prove Theorem~\ref{thm:second}. We begin first by some properties 
of $X_{\alpha}$ and $\varphi$ in the following lemma.
\begin{lem}\label{lem}
	Consider $\eta_{1}(x)\in X_{\alpha}$, $u\in W^{1,p}_{0}(\Omega)$ and 
	a sequence $(u_{n})_{n}\subset W^{1,p}_{0}(\Omega) $ then \\
	\textup{(a)} $\eta_{1}(x)\varphi(u(x))\in L^{1}(\Omega)$.\\
	\textup{(b)} If $u_{n}\rightharpoonup u$ in $W^{1,p}_{0}(\Omega)$, then 
	$\eta_{1}(x)\varphi(u_{n}) \to \eta_{1}(x)\varphi(u)$ in 
	$L^{1}(\Omega)$.   \\
	\textup{(c)} If $||u_{n}||_{1,p} \to +\infty$ and $v_{n}={u_{n}}/{||u_{n}||_{1,p}} 
	\rightharpoonup v$ in $W^{1,p}_{0}(\Omega)$ and almost everywhere 
	in $\Omega$ then 
	${\eta_{1}(x)\varphi(u_{n})}/{||u_{n}||_{1,p}^{p}}\to 0$ in 
	$L^{1}(\Omega)$. 
\end{lem}
\begin{proof}[Proof of Lemma~\ref{lem}]
	Consider only the case $p<N$ (When $p\geq N$, the proof is immediate.)
	
	Consider $q>(p^{*}/\alpha)'$  such that 
	$\eta_{1}(x)\in L^{q}(\Omega)$, there is $\beta>\alpha$ such that 
	$q>(p^{*}/\beta)'>(p^{*}/\alpha)'$.  Set $q_{1}=\beta q'$ and 
	$g\colon \Omega\times\Rset \to\Rset$ defined by 
	$g(x,s)=\eta_{1}(x)\varphi(x)$. Then, $g$ is a Carath\'eodory 
	function and maps $L^{1}(\Omega)$ into $L^{1}(\Omega)$, which 
	follows immediately from the property (iv) of $\varphi$.
	Since $q_{1}<p^{*}$, we have $W^{1,p}_{0}$ embeds compactly in 
	$L^{q_{1}}(\Omega)$, hence we get (a), (b).
	
	The property (c) is a consequence of the properties (i) and (iv) of $\varphi$.
\end{proof}
\begin{proof}[Proof of Theorem~\ref{thm:second}]
	The functional $\Phi$ is well defined and takes its values in $\Rset 
	\cup \{+\infty\}$. Indeed, let $\varepsilon>0$, by ($f_{0}$) and 
	(G$_{2}$), there is $d_{\varepsilon}(x)\in L^{1}(\Omega)$  such that 
	\begin{equation}
		F(x,s) \leq \frac{\lambda_{1}}{p}|s|^{p} + 
		(\eta(x)+\varepsilon)\varphi(s)+d_{\varepsilon}(x)\text{ a. e. in 
		}\Omega \text{ and }\forall s\in\Rset.
		\label{6-1}
	\end{equation}
	For $u\in W^{1,p}_{0}(\Omega)$, we have 
	\begin{equation}
		\Phi(u)\geq \frac{1}{p}\int_{\Omega}\left(|\nabla u|^{p} 
		-\frac{\lambda_{1}}{p} |u|^{p} \right) - \int_{\Omega} \left( 
		\eta_{\varepsilon}(x) \varphi(u(x)) -d_{\varepsilon}(x) \right) 
		-\langle h,u \rangle
		\label{6-2}
	\end{equation}
	where $\eta_{\varepsilon}(x)=\varepsilon+\eta(x)\in X_{\alpha}$.
	By the property (a) of Lemma~\ref{lem}, we get the result.
	
	The functional $\Phi$ is w.l.s.c., it suffices to use \eqref{6-1}, 
	the property (b) of Lemma~\ref{lem} and Fatou's lemma. $\Phi$ is also coercive. 
	Suppose by contradiction that there exists $A\in\Rset$, $(u_{n})_{n} 
	\subset W^{1,p}_{0}(\Omega)$ such that $||u_{n}||_{1,p}\to \infty$ 
	and $\Phi(u_{n}) \leq A$ for all $n\in \Nset$. Set 
	$v_{n}={u_{n}}/{||u_{n}||_{1,p}}$. We can suppose that $v_{n} \rightharpoonup 
	v$ in $W^{1,p}_{0}(\Omega)$, $v_{n}\to v$ in $L^{p}(\Omega)$ and 
	almost everywhere in $\Omega$. 
	
	\textbf{Affirmation} 
	$$
	v=\varphi_{1}\quad\text{or}\quad v=-\varphi_{1}.
	$$
	Indeed, it suffices to divide~\eqref{6-2} by 
	$\frac{1}{p}||u_{n}||_{1,p}^{p}$. Tending $n\to\infty$ and using the 
	property (c) of Lemma~\ref{lem}, we obtain $0\geq 1-\lambda_{1} \int 
	|v|^{p}$, and hence 
	$$
	||v||_{1,p}^{p} \leq 1 \leq \lambda_{1} \int |v|^{p} \leq ||v||_{1,p}^{p}.
	$$
	The affirmation is then a consequence of the variational 
	characterization of $(\lambda_{1},\phi_{1})$. Suppose that 
	$v=-\varphi$ (we get the same conclusion if $v=\varphi_{1}$).
	
	Since $\frac{1}{p}\int |\nabla u_{n}|^{p} - \frac{\lambda_{1}}{p} 
	\int |u_{n}|^{p} \geq 0$, then $A\geq \Phi(u_{n}) \geq -\int 
	G(x,u_{n}) -\langle h,u_{n}\rangle$. Divide by $\varphi(||u_{n}||)$, 
	when $n\to\infty$, since 
	${||u_{n}||}/{\varphi(||u_{n}||)}\to_{n\to\infty}0$ (property 
	(ii) of $\varphi$) when $n$ goes to $\infty$, we get 
	$$
	\liminf_{n} \int_{\Omega} 
	\frac{G(x,||u_{n}||v_{n})}{\varphi(||u_{n}||)}\geq 0.
	$$
	By \eqref{6-1}, the property (iv) of $\varphi$ and Fatou's lemma, we 
	obtain then 
	$$
	0\leq \int_{\Omega} \limsup_{n} \frac{G(x,||u_{n}||v_{n})}{\varphi(||u_{n}||)}.
	$$
	
	Since $\varphi(|s|)=\varphi(s)$, and by the property (iii) of 
	$\varphi$ we have for almost every $x\in\Omega$ that
	$$
	\begin{array}{rl}
		\limsup\limits_{n} \dfrac{G(x,||u_{n}||v_{n})}{\varphi(||u_{n}||)} & 
		=\limsup\limits_{n} \dfrac{G(x,||u_{n}||v_{n}) \varphi(||u_{n}|| 
		v_{n})}{\varphi(||u_{n}|| v_{n}) \varphi(||u_{n}||)}  \\[2mm]
		 &\leq G_{\varphi}^{-}(x).(\varphi_{1}(x))^{\alpha}.
	\end{array}
	$$
	It follows that $0\leq \int_{\Omega} 
	G_{\varphi}^{-}(x).(\varphi_{1}(x))^{\alpha}$, a contradiction with 
	the hypothesis ($G_{2}'$).
\end{proof}
In the following three examples, we can see that neither of the three 
theorems implies the two others.
\begin{ex}
Let $1\leq \alpha\leq p$, $\varphi$ a comparison function of order 
$\alpha$ and $\eta(x)\in X_{\alpha}$. Set 
$F(x,s)=\frac{\lambda_{1}}{p}|s|^{p}+\eta(x)\varphi(x)$. It is clear 
that under the condition 
$\int_{\Omega}\eta(x)(\varphi(x))^{\alpha}<0$, 
Theorem~\ref{thm:second} applies but Theorem~\ref{thm:7} do not 
because $\Gamma^{\pm}_{p}(x)=0$ almost everywhere in $\Omega$. If 
moreover, we suppose $\left|\big\{x\in \Omega;\ 
\eta(x)>0\big\}\right|\neq 0$ then Theorem~\ref{thm:3} do not apply 
neither because the property $(G_{3})$ is not satisfied.
\end{ex}

\begin{ex}
Set $F(x,s)=\frac{\lambda_{1}}{p}|s|^{p}+\eta(x)|s|$ where  $\eta(x)\in 
Y_{1}$. We can check easily that Theorem~\ref{thm:3} applies under 
the condition $(G_{3}')$ which is equivalent to
$$
\int_{\Omega} \eta(x)\varphi_{1}(x) < - |\langle h, \varphi_{1} 
\rangle |.
$$
But the two other theorems do not apply.
\end{ex}

\begin{ex}
Let $a(x)$ a functional in $\mathcal{C}^{\infty}_{c}(\Omega)$ such 
that $a(x)\leq 0$ for all $x\in \Omega$, $\left|\big\{x\in \Omega;\ 
a(x)=0\big\}\right|\neq 0$ and $\left|\big\{x\in \Omega;\ 
a(x)<0\big\}\right|\neq 0$.  Set for some comparison function 
$\varphi$ of order $\alpha$, $1\leq \alpha\leq p$
$$
F(x,s) =\left( \frac{\lambda_{1}}{p}+a(x)\right) + 
|s|^{p}+\big(\varphi(s)|s|^{p}\big)^{1/2}.
$$
Then, $(G_{1})$ and $(G_{1}')$ are satisfied, while $(G_{2})$ and 
(G$_{3}$) are not because 
$G_{\varphi}^{+}=G_{\varphi}^{-}=G_{1}^{+}=G_{1}^{-}=+\infty$ in the 
set of positive measure $\big\{x\in \Omega;\ 
a(x)=0\big\}$.
\end{ex}

\begin{rk}
1) The comparison functions are useless in the autonomous case 
$f(x,s)=f(x)$.  Indeed we have \\
(1.a) $(G_{1})$, $(G_{1}')$ imply $(G_{2})$, $(G_{2}')$ and 
$G_{\varphi}^{\pm}=-\infty$.\\
(1.b) $(G_{2})$, $(G_{2}')$ imply $(G_{3})$, $(G_{3}')$ for any $h\in 
W^{-1,p'}(\Omega)$ and $G_{1}^{\pm}=-\infty$.\\
(1.c) The condition
\begin{equation}\label{G0} 
\lim_{s\to \pm \infty}\frac{G(s)}{|s|}=-\infty
\end{equation}
implies $(G_{3})$, $(G_{3}')$ for any $h\in W^{-1,p'}(\Omega)$ and 
$G_{1}^{\pm}=-\infty$.\\

2) Always, in the autonomous case, we can check easily that the 
condition \eqref{(F)} implies the condition \eqref{G0}.  The converse 
is false.  Indeed, for $f(s)=\lambda_{1}|s|^{p-2}s - \beta |s|^{\beta-2} 
s$ where $1< \beta <p$ so that $F(s)=\lambda_{1}|s|^{p}/p 
-|s|^{\beta}$ and $G(s) = -|s|^{\beta}$.  We have that the function 
$G$ satisfies \eqref{G0} but $\limsup_{s\to \pm 
\infty}{p.F(s)}/{|s|^{p}} = \lambda_{1} $.
\end{rk}

\textbf{Acknowledgement.  }The author is grateful to an anonymous 
referee who brought to his attention some misprints.  This led to a 
better and more correct form of the present paper.

\end{document}